%% file: report.tex
\documentclass[a4paper,11pt,twoside,reqno]{article}
\usepackage[pagebackref,colorlinks=true,urlcolor=blue,linkcolor=blue,citecolor=blue]{hyperref}
\usepackage{a4wide}
\usepackage{times}
\usepackage{amsmath}
\usepackage{amssymb}
\usepackage{amsfonts}
\usepackage{amscd}
\usepackage{amsthm}
\usepackage{amsxtra}

\usepackage[all]{xy}  \CompileMatrices
\usepackage{mathrsfs}
\usepackage{nicefrac}
\usepackage{graphicx}
\usepackage{color}
\usepackage{enumerate}
\usepackage{url}

\input{defs.common}

\newcommand{\R}{\mathbb{R}}

\newcommand{\Q}{\mathbb{Q}}
\DeclareMathOperator{\Tr}{Tr}

\DeclareMathOperator{\frkS}{rk_{\scriptscriptstyle{S}}}

\begin{document}
\title{Dagstuhl Report 13082:\\
  Communication Complexity, Linear Optimization, and lower bounds for the nonnegative rank of matrices}
\author{LeRoy Beasely \and Troy Lee \and Hartmut Klauck \and Dirk Oliver Theis}
\date{Fri, May 17, 2013}
\maketitle

\begin{abstract}
  This report documents the program and the outcomes of Dagstuhl Seminar 13082 ``Communication Complexity, Linear Optimization, and lower bounds for the nonnegative rank of matrices'', held in February 2013 at Dagstuhl Castle.
\end{abstract}

\input{executive-summary}
\input{overview-talks}
\input{problem-sessions}
\input{conclusion}
\input{report.bbl}

\end{document}

%% file: defs.common.tex
\renewcommand{\em}{\sl}

\DeclareMathOperator{\rk}{rk}

\newcommand{\lt}{\left}
\newcommand{\rt}{\right}

\newcommand{\MM}{\mathbb{M}}

\newcommand{\RR}{\mathbb{R}}

\newcommand{\Nm}[1]{{\lt\lVert #1 \rt\rVert}}

\newlength{\algotabbingwidth}


%% file: executive-summary.tex

\section{Executive Summary}

The nonnegative rank is a measure of the complexity of a matrix that has applications ranging from Communication Complexity to Combinatorial Optimization.  At the time of the proposal of the seminar, known lower bounds for the nonnegative rank were either trivial (rank lower bound) or known not to work in many important cases (bounding the nondeterministic communication complexity of the support of the matrix).

Over the past couple of years in Combinatorial Optimization, there has been a surge of interest in lower bounds on the sizes of Linear Programming formulations.  A number of new methods have been developed, for example characterizing nonnegative rank as a variant of randomized communication complexity.  The link between communication complexity and nonnegative rank was also instrumental recently in proving exponential lower bounds on the sizes of extended formulations of the Traveling Salesman polytope, answering a longstanding open problem.

This seminar brought together researchers from Matrix Theory, Combinatorial Optimization, and Communication Complexity to promote the transfer of tools and methods between these fields.  The focus of the seminar was on discussions, open problems and talks surveying the basic tools and techniques from each area.

In the short time since the seminar, its participants have made progress on a number of open problems.


%% file: overview-talks.tex

\section{Overview of Talks}

Background lectures on the connection between matrix factorizations to Communication Complexity and to Combinatorial Optimization were given by the organizers.  More importantly, a number of participants contributed their latest research on factorization ranks.  In this section, we summarize these talks.  The abstracts were provided by the speakers, and have been edited by the organizers.  All errors, inaccuracies, and omissions are due to the editing process.

\subsection{Extended Formulations and Linear Optimization}

\paragraph{Hamza Fawzi}
Many lower bounds on the nonnegative rank only make use of the zero/nonzero pattern of the matrix.  For certain applications, in particular for the extended formulation size lower bounds for approximation problems, nonnegative rank lower bounds need to be shown for matrices that are strictly positive.  Hamza discussed an interesting approach to nonnegative rank lower bounds via conic programming that does not only rely on the zero/nonzero structure of the matrix.  The bound is in many ways analogous to the trace norm lower bound for rank, but making use of the stronger fact that the factorization is nonnegative leads to a copositive program rather than a semidefinite one. For computing the bound in practice, Hamza discussed ways to approximate the bound by semidefinite programs, and examples of using this in practice.

\paragraph{Sam Fiorini}
There is a rich theory on the hardness of approximating NP-optimization problems up to certain factors, given complexity assumptions like P$\neq$ NP. Very recently a similar topic has emerged in the study of polytopes. Sam talked about tradeoffs between the approximation ratio and the size of linear formulations. One notable result in Sam's talk was that approximating CLIQUE to within $n^{1/2-\epsilon}$ requires extended formulations of exponential size. 

\subsection{Complexity}

\paragraph{Nati Linial}
On the first day, Nati Linial treated us to a survey of higher dimensional analogs of familiar combinatorial objects.  For example, we are very familiar with permutation matrices, those matrices with entries from $\{0,1\}$ with exactly one $1$ in every row and column, and know that there are $n! =((1+o(1))n/e)^n$ many of them.  What about $3$-dimensional tensors with entries from $\{0,1\}$ and exactly one $1$ along every row, column and shaft?  Such $2$-dimensional permutations turn out to coincide with latin squares and it is known that there are $((1+o(1))n/e^2)^{n^2}$ many of them.  This relies on some beautiful work on the minimum permanent of doubly stochastic matrices.  Nati conjectures that the formula generalizes to count the number of $d$-dimensional permutations, described by a $d+1$-tensor with one $1$ along every line.  That is, that the number of $d$-dimensional permutations is $((1+o(1))n/e^d)^{n^d}$.  He is able to show such an upper bound, but the lower bound remains open.

\paragraph{Sebastian Pokutta}
In order to prove that extended formulations for approximating optimization problems need to be large, communication and information complexity are important tools.  In his talk Sebastian described a new approach on how to prove lower bounds on the nonnegative rank of matrices corresponding to the unique disjointness problem when perturbed.  He gave tight lower bounds using a new information theoretic fooling set method.

Since the seminar, Sebastian and his co-author Gab\'or Braun have made available a preprint containing these results \cite{BraunPokutta13}.

\paragraph{Hans Raj Tiwary}
There are entire books of NP-complete problems and explicit reductions between them.  For the extension complexity of the associated polytopes, however, this book is still slowly being written---usually by arguing that $P$ is a projection of $Q$ or finding $P$ as a face of $Q$.  Hans discussed the intriguing possibility of automatically turning an NP gadget reduction into a polytope reduction.  While still not a general theory, Hans can currently do this for many NP-hard problems and their associated polytopes.

\paragraph{Nicolas Gillis}
Nicolas spoke about the problem of actually computing a non-negative factorization of a nonnegative matrix.  This talk was important to seminar participants on small matrices, allowing them to test the quality of their lower bounds against upper bounds.  On small matrices, these upper bounds can be found computationally.  The problem also has applications to compression of images, to identifying topics in documents, even to identifying the mineral composition of rocks from spectral data (hyper-spectral imaging).  Nicolas discussed specifically the case of separable matrices.  An $n$-by-$n$ matrix $M$ is $r$-separable if it has a factorization $M=WH$ where $W$ is $n$-by-$r$, $H$ is $r$-by-$N$ and moreover $W$ is a subset of the columns of $M$.  Such types of factorization can be more useful in practice.  Nicolas talked about a linear programming approach to this problem that is polynomial time and moreover outperforms previous approaches in practice.

\subsection{Matrix Theory}

\paragraph{Alexander Guterman}
Alexander Guterman gave a survey talk on various matrix ranks over semirings.  A big focus was on tropical algebra over the real number with operations $a \oplus b=\max{a,b}$ and $a \otimes b=a+b$.  Tropical algebra provides a way of formulating many hard combinatorial optimization problems (like scheduling problems) in terms of a very elegant linear algebraic type language.  In tropical linear algebra there are varying notions of linear independence, for example Gondran-Minoux independence, weak linear independence, and strong linear independence.  Each of these gives rise to a different notion of rank of a matrix and a hierarchy of these ranks is known.

\paragraph{Yaroslav Shitov}
Yaroslav continued talking about tropical matrix rank, in particular the tropical factorization rank.  This is defined as the minimum $k$ such that $A=B \bigotimes C$ for a $n$-by-$k$ matrix $B$ and $k$-by-$n$ matrix $C$. Note that in tropical matrix multiplication $(B \bigotimes C)(i,j)=\min_t B(i,t)+C(t,j)$.  Yaroslav mentioned a very interesting application of the tropical factorization rank.  Say that we are given an instance of the traveling salesman problem, with distances specified by a matrix $A$, and moreover we are given a tropical factorization $A=B \otimes C$ that witnesses that $A$ has constant factorization rank.  Then the resulting traveling salesman instance can be solved in polynomial time!  This is a result of Barvinok, Johnson, Woeginger, and Woodroofe.  Yaroslav also showed that the problem of detecting if the tropical factorization rank of a matrix is at most $8$ is NP-hard.

\paragraph{Richard Robinson}

In his talk, Richard Robinson gave a characterization, among all nonnegative matrices, of the extreme-ray / facet slack matrices of polyhedral cones, and vertex/facet slack matrices of polytopes.  This characterization leads to an algorithm for deciding whether a given matrix is a vertex/facet slack matrix.  The underlying decision problem is equivalent to the polyhedral verification problem whose complexity is unknown.


%% file: problem-sessions.tex

\section{Problem discussion sessions, and subsequent developments}

Here we report on the status of questions which were presented during the Problem Sessions.  Near the end of this section, we discuss developments on problems, which were not presented during the Problem Sessions, but discussed during the seminar.

\subsubsection{Real vs.\ rational nonnegative rank.}\label{prob:real-vs-rational} 
Presented by Dirk Oliver Theis; problem based on a problem by Cohen and Rothblum from 1991.  Give a non-trivial bound for $\rk_{\Q_+}(A) - \rk_{\R_+}(A)$!  For example, is it true that for every rational nonnegative matrix~$A$ we have $\rk_{\Q_+}(A) \le \rk_{\R_+}(A) + 1$?

The original question asks for equality between the two ranks, but currently no non-trivial bounds for $\rk_{\Q_+}(A) - \rk_{\R_+}(A)$ or even $\rk_{\Q_+}(A) \bigm/ \rk_{\R_+}(A)$ are known.

In the discussion, Nati Linial pointed to Micha Perles discovery of non-rational polytopes, and the studies by Richter-Gebert and others of the realization spaces of polytopes.

\medskip%
A related problem is the following.
\medskip%

\noindent
\textbf{Complex vs.\ real positive semidefinite rank.}  One can ask a similar question for the positive semidefinite rank.  The positive semidefinite rank over $\mathbb{R}$ of a matrix $A \in \mathbb{R}^{m \times n}$ is the minimal $r$ such that there are $B_i \in \mathbb{R}^{r \times r}$ for $i=1,\ldots, m$ and $C_j \in \mathbb{R}^{r \times r}$ for $j=1, \ldots, n$ such that $A(i,j)=\Tr(B_i^*C_j)$.  The positive semidefinite rank over $\mathbb{C}$ is defined analogously with $B_i, C_j \in \mathbb{C}^{r \times r}$.  Is the positive semidefinite rank over $\mathbb{R}$ equal to positive semidefinite rank over $\mathbb{C}$?  This is the simplest in a family of questions: The same needs to be asked for rational vs.\ real positive semidefinite rank.  In the discussion, Nicolas Gillis pointed out that this question also must be settled for the copositive ranks.  In this case one looks for a factorization $A(i,j)=\Tr(B_i^* C_j)$ where each $B_i$ is copositive and each $C_j$ is completely positive.

\subsubsection{Square-root rank.}\label{prob:squarerootrank}
Proposed by Richard Robinson.  Given nonnegative matrix $M \in \R^{p \times q}_+$, we say that $A \in \R^{p \times q}$ is a Hadamard square-root if $(A)^2_{ij} = M_{ij}$.  Define $\rk_{\text{\tiny$\sqrt\cdot$}}(M)$ as the minimum rank of~$A$ such that~$A$ is Hadamard square-root of~$M$.  This is equivalent to a version of the positive semidefinite rank where the matrices in the factorization are constrained to have rank~$1$.  Hence, $\rk_{PSD} (M) \leq \rk_{\text{\tiny$\sqrt\cdot$}}(M)$ holds.  For vertex-facet slack matrices of polytopes, $\rk(S) \leq \rk_{\text{\tiny$\sqrt\cdot$}} (S)$.  A large number of observaions have led Richard to conjecture the following.
\begin{quote}
  If $\rk(S) = \rk_{\text{\tiny$\sqrt\cdot$}}(S)$ holds for the slack matrix~$S$ of a polytope, then the entries in the hadamard square-roots in $\rk_{\text{\tiny$\sqrt\cdot$}}(S)$ can be taken to be nonnegative.
\end{quote}
In the discussion, Samuel Fiorini suggested to look specifically at the matrix $M_{ab} = (1-a^Tb)^2$.

\subsubsection{Positive semidefinite rank of matrices defined by polynomials.}\label{prob:psdrkpolynomial}
Proposed by Troy Lee.  What is the positive semidefinite rank of the matrix
\begin{equation*}
  M(x,y) = (x^ty-1)(x^ty-2),
\end{equation*}
where $x, y$ range over all $\{0,1\}^n$?.  The motivation is that such a matrix~$M$ is a submatrix of the slack matrix of the correlation polytope.  One can define a whole family of submatrices of the slack matrix of the correlation polytope by taking a quadratic polynomial~$p$ which is nonnegative on nonnegative integers, and letting $M(x,y) = p(\lvert x \cap y\rvert)$.  In the discussion, Sam Fiorini pointed out that to show a strong lower bound on $M$ one would have to focus on more than just the entries of the matrix which take values in some small set.  He also mentioned that this matrix can be approximated by one that does have low positive semidefinite rank, namely the matrix $N(x,y)=(x^ty-3/2)^2$.

A toy version of this problem asks about the positive semidefinite rank of the $n$-by-$n$ matrix $M_n(i,j)=(i-j-1)(i-j-2)$.  Seminar participant Jo\~{a}o Gouveia \cite{Gouveia13} answered a question of Lee and Theis \cite{LeeTheis12} by showing that the psd rank of $M_n$ goes to infinity with $n$.

\subsubsection{A query complexity problem.}\label{prob:querycomplexity}
Proposed by Raghav Kulkarni.  For $f\colon \{0,1\}^n \to \{0,1\}$ and $z\in\{0,1\}^n$, we say that $i$th bit of~$z$ is \textit{sensitive} if $f(z_1,\ldots,\bar{z_i},\ldots,z_n) \neq f(z_1,\ldots, z_i,\ldots,z_n)$.  Let $s(f,z)$ be the number of sensitive bits of~$z$ and $s(f) = \max_z \{s(f,z)\}$ the maximum number of sensitive bits of any argument.  These concepts arise in the context of decision tree complexity.  For $x, y \in \{0,1\}^n$ let $f(x,y)$ be a 2-parameter function, and let $M_f(x,y)=f(x,y)$ be the corresponding matrix.  Raghav conjectures that 
$\log \rk_\R(M_f) \leq \mathrm{poly}(s(f))$.

In the following discussion, Hartmut Klauck asked about block sensitivity and Raghav said the conjecture is true with sensitivity replaced by block sensitivity.  Hartmut also suggested easier versions of the conjecture where, for example, the rank is replaced by sign rank which is the minimum rank of a matrix that entrywise agrees with the target matrix in sign.  Nati Linial then asked if assuming the log-rank conjecture is true implies anything for this conjecture.

\subsubsection{A rigidity-type question.}\label{prob:adi}
Proposed by Adi Shraibman.  The famous matrix rigidity problem of Valiant asks to explicitly construct matrices with high rank and such that if a constant fraction of the entries are arbitrarily changed, the rank remains high. While probabilistic constructions exist, finding explicit constructions remains a hard open problem.

Consider the following variant of the problem known as discrepancy games.  Here you start with an empty $n \times n$ matrix.  Two players, Balancer ($+1$) and Unbalancer ($-1$), take turns assigning entries of the matrix to their associated value.  Balancer wants to make all combinatorial rectangles balanced, while Unbalancer wants to make them unbalanced. In this game it is known that Balancer can get ensure an upper bound of $s^{3/4}$ on discrepancy after $s$ rounds.

Here is another variant that is open.  In this case we begin with a $\{-1,+1\}$ valued matrix with discrepancy $n^{3/2}$.  Say a Hadamard matrix. Balancer picks certain $+1$'s. Unbalancer picks certain $-1$'s. Over the course of the game, can Balancer maintain discrepancy to be less than $s^{3/4}$ after $s$ moves?  In this variation you cannot rely on strategy stealing.

\subsubsection{Extension complexity of stable set polytopes of split-graph-free perfect graphs.}\label{prob:stable-set}
Proposed by Samuel Fiorini.  A split graph is a graph in which the vertices can be partitioned into a clique and an independent set.  Let~$H$ be a split graph, and consider the class of all graphs~$G$ not containing~$H$ as an induced subgraph.  What is the extension complexity of the stable set polytopes of this class of graphs?  This problem is motivated by a recent result of Bousquet, Lagoutte, and Thomass\'{e} \cite{BLT13}.  They provide a certificate of size $O(\log(n))$ proving that a clique and an independent set do not intersect for $H$-free graphs.


\subsubsection{Matrices with low non-negative rank are a low-dimensional subset in the manifold of rank-bounded matrices.}\label{prob:dimsemialg} 
It is an easy fact that the nonnegative rank is semicontinuous: If $A_j$ tends to $A$, then $\rk_+(A) \leq \liminf \rk_+(A)$.  Let $n,k$ be nonnegative integers such that $3 \le k \ll n$, and consider the set of $n\times n$ matrices
\begin{equation*}
  \{A \in \mathbb{R}_+^{n \times n} \mid \rk A = k\text{, } \rk_+ A \le n-1\}.
\end{equation*}
Does this set contain interior points within the manifold $\{A \in \mathbb{R}_+^{m \times n} \mid \rk A = k\}$ of nonnegative rank-$k$ matrices?  (The ``$n-1$'' is somewhat arbitrary, and should be replaced by an appropriate function of~$n$.)  This question asks for the dimension of the set of matrices of ``small'' nonnegative rank as a semialgebraic subset of the variety of rank-$k$ $n\times n$ matrices.  In the discussions, people expressed that the intuitively obvious answer to the question is yes.  But a recent result of Yaroslav Shitov proves that, for $k=3$, \textit{every} such matrix~$A$ has nonnegative rank at most $6n/7$ \cite{Shitov13}.  However, the question with ``$n-1$'' may still be true for large~$k$.

\subsubsection{Fooling-sets and rank.}\label{prob:fooling}  
Let~$A$ be an $n \times n$ matrix over a field~$\mathbb K$ satisfying $A_{kk} = 1$ for all $k$ and $A_{k\ell}A_{\ell k}=0$ whenever $k \neq \ell$.  Dietzfelbinger et al.\ (1996) proved that $n \leq \rk(A)^2$.  The question raised by Diezfelbinger et al.\ is whether this bound is asymptotically ($n \to \infty$) tight.

This problem was fully settled by Friesen and Theis in the case of nonzero characteristic shortly before the seminar 
\cite{FriesenTheis13}.  In summary, the following is known.
\begin{center}
  \begin{tabular}{r|ccc}
    ~        & \multicolumn{3}{l}{characteristic of $\mathbb K$}\\
    ~        & $0$ & $2$ & $\ge 3$\\
    \hline
    $A$ 0/1 entries       & open & tight & open\\
    $A$ arbitrary entries & open & tight & tight
\end{tabular}
\end{center}
The major open question is when the characteristic is~0 and the entries are arbitrary.  At the time of the seminar, the best separation was by Klauck and de Wolf \cite{KlauckWolf13} who gave an example where $\rk(A) \le n^{0.613\ldots}$ with integral entries of small modulus.  After the seminar, in the case of characteristic zero, Troy Lee was able to improve the best known bound to $\rk(A) \le n^{0.594\ldots}$ with a method that warrants future investigation.

\subsubsection{Polygons --- or, more generally, rank-3 matrices.}

A problem which was discussed intensely during the seminar was the extension complexity of polygons, where, at the time of the seminar, a lower bound of $\Omega(\sqrt n)$ was known, and the trivial upper bound~$n$.

The following problem by Beasley \& Laffey \cite{BeasleyLaffey09} is both more general and more specific: Given any sub-semiring $S$ of~$\RR_+$ and $n\ge 6$, is there a matrix $A \in \MM_{n\times n}(S)$ such that $\rk A = 3$ and $\frkS(A) = n$?  At the time of the seminar, this was open even for~$S=\RR_+$.  Following the seminar, participant Yaroslav Shitov \cite{Shitov13} has settled this problem for the semiring $\R_+$: every rank-3 nonnegative $n\times n$ matrix with $n\ge 7$ has nonnegative rank at most $6n/7$.

\subsubsection{Euclidean distance matrices.}

Do ``generic''/``random'' Euclidean distance matrices\footnote{A $d$-dimensional, Euclidean distance matrix of size~$n$ is defined by points $x_1,\dots,x_n$ in $d$-dimensional Euclidean space.  Its entries are $\bigl( \Nm{ x_k - x_\ell} )_{k,\ell}$.} have full nonnegative rank?\footnote{An affirmative proof of this for $d=1$ by Lin and Chu (2011) is fatally flawed.}  

For $d=1$, Shitov's above mentioned result gives a negative answer to the question.

\subsubsection{How does nonnegative rank behave under tensor products?}
This question is interesting on its own, and also has application to communication complexity.  A very strong conjecture, discussed at the workshop, would be that the nonnegative rank is multiplicative, as the rank is.  That is, that $\rk_+(A \otimes B)=\rk_+(A)\rk_+(B)$.

Collaboration between two seminar participants, Nicolas Gillis and Hamza Fawzi showed that this strong conjecture is false.  Specifically, Nicolas wrote software for computing the nonnegative rank of small-size matrices, and using this software Hamza was able to disprove the conjecture.

\subsubsection{Vertex/facet slack matrices vs.\ general nonnegative matrices.}

Some lower bounds can be proved for matrices which arise from vertex/facet slack matrices of polytopes.  It is an open question whether some of the bounds behave fundamentally different in the case of vertex/facet slack matrices (cf.\ e.g., \ref{prob:squarerootrank}).  As a first step towards resolving this type of questions, a linear-algebraic characterization of these matrices was obtained by Jo\~ao Gouveia, Richard Robinson, and Rekha Thomas, and presented during the seminar.  It turned out that Volker Kaibel, Roland Grappe, and Kanstantsin Pashkovich had a similar approach.  Their results were combined in the recent paper \cite{GouveiaGrappeKaibelPashkovichRobinsonRhomas14}.

\subsubsection{Extensions/factorizations over the positive semidefinite cone}

In the year preceding the seminar, matrix factorizations over other cones than $(\RR_+)^r$ have gained importance, specifically the cone of positive semidefinite matrices.  From the combinatorial optimization point of view, a very basic question there is, whether there exist polytopes whose extension complexity is exponential in the dimension --- the same question for $(\RR_+)^r$ was settled by Rothvo\ss~\cite{Rothvoss12}.  In a very recent paper, seminar participant Sebastian Pokutta, together with two coauthors, Jop Bri\"et and Daniel Dadush, have answered that question in the affirmative \cite{BrietDadushPokutta13}.

In the context of factorizations over other cones, recently, the conference participants Samuel Fiorini and Hans Raj Tiwary~\cite{FioriniTiwary13} have observed that a 2012 theorem by Alexander Maksimenko~\cite{Maksimenko12} together with a result by Samuel Burer~\cite{burer2009} implies that every 0/1 polytope whose vertex set can be described by a polynomial predicate has a polynomial sized copositive extension.


%% file: conclusion.tex

\section{Conclusion}

A natural approach to solving hard combinatorial optimization problem is to give a formulation as a linear program and solve it using standard techniques.  An important topic initiated by Yannakakis is to investigate the size of extended formulations of optimization problems.

Recently the theory of representing hard optimization problems via extended formulations has seen much progress.  Techniques from communication complexity and matrix theory have been essential to investigate how large extended formulations need to be, finally improving on Yannakakis' seminal results.

The seminar brought together researchers from the areas of optimization theory, complexity theory, and matrix theory, to further collaboration on these and newly emerging topics.  Exciting progress was reported on proving lower bounds for the nonnegative rank, on the hardness of approximation using extended formulations, and on new notions of matrix ranks.

For the future we hope that similar progress will soon be made on the topic of using semidefinite programming to solve hard optimization problems.  Intriguingly this problem is connected to quantum communication complexity.

Dagstuhl provided a wonderful environment for many informal discussions as well as talks, plus an exciting open problems session.  The opportunity to have this seminar was well appreciated by the participants, many of them who were new to the center.


%% file: report.bbl
\newcommand{\etalchar}[1]{$^{#1}$}